%% file: article5a.tex
\documentclass[12pt,leqno]{article}

\usepackage{amssymb}
\usepackage{amsmath}
\usepackage{graphicx,psfrag,color}
\usepackage{url}
\usepackage{pstool}

\begin{document}

\title{Nonrepresentable relation algebras from group
systems}
\author{Andr\'eka, H., Givant, S. and N\'emeti, I.}%
\date{}
\maketitle

\input article5a-def.tex
%\keywords{relation algebra, group, coset, measurable atom, Boolean
%algebra}
%\subjclass[2010]{Primary: 03G15; Secondary: 20A15}
\begin{abstract} A series of nonrepresentable relation algebras is
constructed from groups. We use them to prove that there are
continuum many subvarieties between the variety of representable
relation algebras and the variety of coset relation algebras. We
present our main construction in terms of polygroupoids.
\end{abstract}

\section{Introduction}\label{intro-sec}

In 1941, Tarski \cite{tarski41}  defined the class $\RA$ of relation
algebras. This was the culmination of previous extensive work by
Peirce and Schr\"oder, who worked towards extending Boole's system
to capture more of our reasoning than propositional calculus does.
In the early 1940s, Tarski proved that set theory can be built on
the equational theory of relation algebras \cite{tg}. It may be just
an accidental fact that mathematics today is based on first-order
logic and not on the equational theory of relation algebras
\cite{morigin, prhist}. Still today, $\RA$ is an active field of
research. It is a central topic in algebraic logic \cite{giv18a,
giv18b, hh02, ma06}, but relation algebras also have deep
connections with other areas of mathematics, such as geometry, graph
theory, combinatorics, group theory, linguistics, arrow logic, modal
logic, to list some \cite{vB, comer0, Du, khal, lyndon, maddux, mv}.
Relation algebras are extensively used in computer science
\cite{bks, prda, sch}.
%jipsenre, Orlowskara, Szabolcs is hivatkozni

In the 1940s, J. C. C. McKinsey noticed that the complex algebra
(the natural algebra of complexes or subsets) of a group is a
relation algebra. J\'onsson and Tarski proved that the complex
algebra of a Brandt groupoid is a relation algebra \cite[section
5]{jt52}. Comer \cite{comer1} defined polygroupoids (these turned
out the same as regular, reversible-in-itself multigroups %with an absolute unit
in Dresher-Ore \cite{dresore}) and proved that relation algebras are
exactly subalgebras of complex algebras of polygroupoids. He showed
that some theorems of relation algebra theory are easier to prove
for polygroupoids first and then transfer the result to relation
algebras \cite{comer2}. He suggested a research project of doing a
portion of relation algebra theory in terms of polygroupoids. In
this paper, we take up his initiative since we find that a large
part of our work can be more conveniently formulated in terms of
polygroupoids than in terms of relation algebra theory. In
particular, we can introduce the class $\RRA$ of representable
relation algebras by generalizing the notion of the Cayley
representation of a group to polygroupoids. Representable relation
algebras are the primary or generic examples of relation algebras.
%dresher-Ore ellenorizni valamikor
%irodalom beszurni: comer.

A generalization of group relation algebras is introduced in
\cite{ga02, giv1}, by using an intricate system of groups and their
factor groups (see Figure 1). They are all representable, their
representations closely follow the Cayley representations of the
groups involved. However, with a slight ``shift" we can make some of
these group relation algebras nonrepresentable. One such example is
presented in \cite[section 5]{andgiv1}. In this paper, we generalize
that example to construct an infinite series of nonrepresentable
shifted group relation algebras (coset relation algebras, in short).
We use these algebras to prove that there are uncountably many
varieties of coset relation algebras that all contain the variety of
representable relation algebras (Theorem \ref{var-thm}).

It is known to be difficult to construct nonrepresentable relation
algebras, because the finitely axiomatized class of relation
algebras approximates surprisingly well the nonfinitely
axiomatizable class of representable relation algebras. Our method
for constructing nonrepresentable relation algebras has novel
features, it differs in intuition from the ones available in the
relation algebra literature. One novel feature is that they are
constructed from systems of groups with only slightly distorting the
group system structure. Only how the groups are connected is
distorted, the groups themselves are kept intact. Another feature is
that the cause of nonrepresentability in these algebras is ``sheer"
structure, not size, not that some parts of the algebras are too big
or too small.

Why does one want to construct nonrepresentable relation algebras,
when the representable ones are of the primary interest? For
example, because they can be used to prove various properties of
representable relation algebras. For example Monk \cite{Monk} uses a
series of nonrepresentable relation algebras to prove that the
variety \RRA\ cannot be axiomatized by a finite set of equations.
Similarly, J\'onsson \cite{J91} proves that each equational axiom
system of \RRA\ must use infinitely many variables. His proof relies
on nonrepresentable relation algebras constructed from projective
planes. For a sample of constructions of nonrepresentable algebras
see \cite{amnsplit, hh01, jon59, ma18, McK, SDis}.

In section \ref{groupoid-sec}, we introduce previous work that we
need from relation algebra theory, in terms of polygroupoids. We
assume basic knowledge of groups, e.g., normal subgroup, factor
group. We hope that this section can be read without any background
in relation algebra. We also hope that some knowledge acquired in
relation algebra theory can be made available this way to a larger
community of mathematicians. In sections \ref{frame-sec} and
\ref{shift-sec} we still work solely with polygroupoids and groups.
In section \ref{frame-sec}, we construct a group frame starting from
any commutative group. In section \ref{shift-sec}, we add the
``shift" and we prove that the algebras obtained are all
nonrepresentable. This is one of the two main theorems of the paper.
Relation algebras are introduced in section \ref{var-sec}, where we
begin to use relation algebra terminology. We show that the examples
we constructed in sections \ref{frame-sec}, \ref{shift-sec} are
diverse enough in that none of them can be embedded into the other.
This feature is used in proving our second main theorem, Theorem
\ref{var-thm}.

\section{Groups, groupoids, group systems}\label{groupoid-sec}

We recall the definition of polygroupoids. The notion of a
polygroupoid is a generalization of that of a group. Brandt
groupoids are generalizations of groups in that the binary
composition operation is partial in them and that they can have more
than one identity elements. Composition can be multivalued in
polygroups, in addition.

Some notation that will be convenient (and is customary) to use in
partial multivalued algebras. When $a\scir b$ is not defined, we
write $a\scir b=\emptyset$. When $X,Y\subseteq M$  we write $X\scir
Y=\{c : c\in a\scir b\ \mbox{ for some }a\in X\ \mbox{ and }b\in
Y\}$, and we often write just $a$ in place of the singleton $\{
a\}$.

\begin{Def}[Polygroupoid, {\cite[Definition 3.1]{comer1}}]\label{poly-def}
By a polygroupoid we understand a structure $\Mm=(M,\scir, I,
{}^{-1})$  where $\scir$ is a partial multivalued binary operation
on $M$, $I\subseteq M$ and ${}^{-1}$ is a unary operation on $M$
that satisfy the following three conditions for all $a,b,c\in M$.
\begin{description}
\item[(i)] $(a\scir b)\scir c=a\scir (b\scir c)$,
\item[(ii)] $a\scir I=a=I\scir a$,
\item[(iii)] $a\in b\scir c$ holds if and only if $b\in a\scir c^{-1}$ and this holds if and only if
$c\in b^{-1}\scir a$.\qed
\end{description}
\end{Def}

Polygroupoids appear in many parts of mathematics. In particular,
groups and Brandt groupoids are special polygroupoids. The latter
can be characterized as those polygroupoids in which the composition
operation $\scir$ is not multivalued.

An element $a$ of a polygroupoid is called a \emph{loop} if there is
$x\in I$ such that $x\scir a\scir x=a$ (intuitively, if the domain
and range of $a$ coincide). A polygroupoid is said to be
\emph{locally functional} if the product of two loops $f,g$ in $\Mm$
cannot have more than one value, i.e., if $|f\scir g|\le 1$. Let
$\PG$ denote the class of all polygroupoids and let $\LPG$ denote
the class of all locally functional polygroupoids.

Surprisingly, locally functional polygroupoids have a rather rich
structure. They are put together from various groups and their
factor groups in an orderly manner. The novelty here is the
appearance of factor groups, they appear because of the multivalued
nature of the composition operation. The notion of a group coset
frame emerges as the complete description of the structure of
locally functional polygroupoids, see Theorem \ref{lpg-thm} below.

For introducing group coset frames, we need some notation concerning
groups. When $\Gg$ is a group, we usually denote its universe by
$G$, its binary operation by $\scir$, its identity element by $e$
and inverse of $g\in G$ by $g^{-1}$. When $H$ is a normal subgroup
of $\Gg$, we denote the set of cosets of $H$ in $\Gg$ by $G\slash
H$. Thus $G\slash H=\{ g\scir H : g\in G\} = \{ H\scir g : g\in
G\}$, this is the universe of the factor group $\Gg\slash H$. More
generally, when $X\subseteq G$ we denote $X\slash H = \{ g\scir H :
g\in X\}$ and $g\slash H = g\scir H$. Since groups are special
polygroupoids, we use the notation introduced for polygropoids for
groups, too. In particular, when $H,K$ are two normal subgroups of
$\Gg$, their complex product, $H\scir K$, is the normal subgroup
generated by them. Now, $H\scir K$ is both a union of cosets of $H$
and a union of cosets of $K$. Often, we will identify $H\scir K$
with $H\slash K$ or $K\slash H$ while the exact connection is
$H\scir K=\bigcup H\slash K$. We hope this will not lead to
confusion.

\begin{Def}[Group coset system, \cite{ga02}, {\cite[section 3]{andgiv1}}]\label{gcs-def}
Let $I$ be any set and $\E\subseteq I\times I$ be an equivalence
relation on $I$. Assume that  $\G=\langle \Gg_x : x\in I\rangle$ is
a system of groups, $\varphi=\langle\varphi_{xy} :
(x,y)\in\E\rangle$ is a system of isomorphisms between their factor
groups,
\[ \varphi_{xy} :\Gg_x\slash H_{xy}\to \Gg_y\slash K_{xy}\quad\mbox{ is
an isomorphism }\] and $C=\langle C_{xyz} : (x,y),(y,z)\in\E\rangle$
is a system of cosets of these factors, in more detail,
\[ C_{xyz}\quad\mbox{ is a coset of }H_{xy}\scir H_{xz} .\]
Then we call $(\G,\varphi,C)$ a \emph{group coset system}.\qed
\end{Def}

Group coset frames are group coset systems in which there is a
concert between the ingredients of the group coset system, see
Figure 1.  When $f$ is a function with domain $A$ and $X\subseteq
A$, we denote the image of $X$ under $f$ by $f[X]=\{ f(x) : x\in
X\}$. When $C$ is a coset of a normal subgroup $H$, we denote the
involution induced by it, a special automorphism, by $\tau(C)$, thus
\[\tau(C)(D)=C^{-1}\scir D\scir C\] for all cosets $D$ of $H$. We will
denote compositions of functions with a sign different from the
usual $\scir$, moreover we will use relation-type composition as in
category theory where the order of applying the functions is
different. In more detail: the composition of binary relations $R,S$
is
\[ R\mid S=\{(u,v) : (u,w)\in R\mbox{ and }(w,v)\in S\mbox{ for some
}w\}.\] When $f:U\to V$ is a function, we consider it as a binary
relation $\{(u,f(u)) : u\in U\}$, then $f\mid g:U\to W$ if $g:V\to
W$ and
\[ (f\mid g)(u) = g(f(u)) .\]

\begin{Def}[Group coset frame, {\cite[Definitions 3.2,
4.1]{andgiv1}}]\label{frame-def}\label{cosetframe-def}  A group
coset system $(\G,\varphi,C)$ is called a \emph{group coset frame}
if the following eight conditions are satisfied for all
$(x,y),(y,z),(z,w)\in\E$.
\begin{description}
\item[(i)]
$\varphi_{xx}$ is the identity function on $G_x\slash \{e_x\}$,
where $e_x$ is the identity element of $G_x$.
\item[(ii)] $\varphi_{yx}$ is the inverse of $\varphi_{xy}$. In
particular, $K_{xy}=H_{yx}$.
\item[(iii)]
$\varphi_{xy}[H_{xz}\slash H_{xy}]=H_{yz}\slash H_{yx}$.
\end{description} Assume that
{(iii)} holds. Define
 $\varphi_{xy}^{z}(g\slash (H_{xy}\scir H_{xz})) =
\varphi_{xy}(g\slash H_{xy})\scir H_{yz}$.
\begin{description}
\item[(iv)]
$\varphi_{xy}^{z}\mid\varphi_{yz}^{x}=\tau(C_{xyz})\mid\varphi_{xz}^{y}$.
\end{description}
We now list the connections between the cosets $C_{xyz}$ that are
required to hold in a group coset frame.
\begin{description}
\item[(v)] $C_{xyy}=H_{xy}$.
\item[(vi)] $\varphi_{xz}[C_{xyz}]=C_{zyx}^{-1}$.
\item[(vii)] $\varphi_{xy}[C_{xyz}]=C_{zyx}^{-1}$.
\item[(viii)] $C_{xyz}\scir C_{xzw}=\varphi_{yx}[C_{yzw}\scir H_{yx}]\scir
C_{xyw}$.\qed
\end{description}
\end{Def}

\begin{figure}\label{frame-fig}
\centering \small \psfrag{C}[c][c]{$C_{xyz}$}
\psfrag{xy}[b][b]{$\Gg_x\slash H_{xy}$}
\psfrag{yx}[b][b]{$\Gg_y\slash H_{yx}$}
\psfrag{yxz}[t][t]{$\Gg_y\slash(H_{yx}\scir H_{yz})$}
\psfrag{x}[r][r]{$\Gg_x$} \psfrag{y}[l][l]{$\Gg_y$}
\psfrag{yz}[l][l]{$\Gg_y\slash H_{yz}$}
\psfrag{zy}[l][l]{$\Gg_z\slash H_{zy}$} \psfrag{z}[t][t]{$\Gg_z$}
\psfrag{f}[b][b]{$\varphi_{xy}$} \psfrag{n}[b][b]{$\varphi^z_{xy}$}
\psfrag{xyz}[t][t]{$\Gg_x\slash(H_{xy}\scir H_{xz})$}
\includegraphics[keepaspectratio, width=0.8\textwidth]{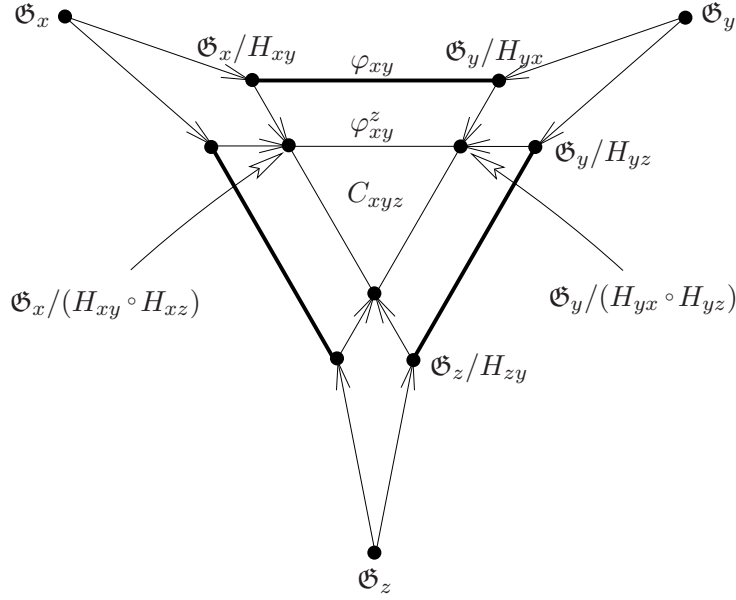}
\caption{Group coset frame $(\langle\Gg_x : x\in
I\rangle,\langle\varphi_{xy} :(x,y)\in\E\rangle,\langle C_{xyz} :
(x,y),(y,z)\in\E\rangle)$.}
\end{figure}

Theorem \ref{lpg-thm} below states that group coset frames
completely capture the structure of locally functional
polygroupoids, in an analogous manner as a system of groups captures
the structure of a Brandt groupoid. Let us call a polygroupoid
\emph{connected} when for all identity elements $x,y$ there is an
element $a$ such that $x\scir a\scir y=a$. It is known that the
structure of each connected Brandt groupoid can be described by a
group $\Gg$ and a set $I$, as follows. To each group $\Gg$ and set
$I$, a Brandt groupoid $\Bb(\Gg,I)$ is defined, such that the
universe is $\{ (x,g,y) : x,y\in I\mbox{ and }g\in G\}$ and
multiplication is $(x,g,y)\scir (y,h,z) = (x,g\scir h, z)$. One can
prove that each connected Brandt groupoid is isomorphic to one of
these $\Bb(\Gg,I)$. When $\Bb$ is not connected, we have an
equivalence relation $\E$ on $I$ and distinct blocks of $\E$ may
have distinct groups.
Now, group coset frames capture the structures of locally functional
polygroupoids, in an analogous way. First we introduce the
polygroupoid analogues of $\Bb(\Gg,I)$.

\begin{Def}[Structure associated to group coset
system]\label{str-def}  As\-sume that $(\G,\varphi,C)$ is a group
coset system. We define the structure $\Pp(\G,\varphi,C)$ as
follows. The universe is
\begin{description}\item{} $P=\{ (x,g,y) :
(x,y)\in\E\mbox{ and } g\in G_x\slash H_{xy} \}$.
\end{description}
Multiplication on $P$ is a multivalued binary partial function
defined as follows. We will denote the multiple $(x,g,y)\scir
(y,h,z)$ simply by juxtaposition $(x,g,y)(y,h,z)$. Assume that
$(x,g,y), (w,h,z)\in P$. Their multiple is defined exactly when
$y=w$ and
\begin{description}
\item{}
$(x,g,y)(y,h,z)=\{ (x, k, z) : k\subseteq
\varphi_{yx}(\varphi_{xy}(g)\scir h)\scir C_{xyz}\ \} $.
\end{description}
The set $E$ of identities of $\Pp$ is
\begin{description}
\item{}
$E = \{ (x,H_{xx},x) : x\in I \} $
\end{description}
and inverse is defined as
\begin{description}
\item{}
$(x,g,y)^{-1}=(y,\varphi_{xy}(g^{-1}),x) $.\qed
\end{description}
\end{Def}

The following theorem states that the concrete polygroupoids defined
in the above definition exhaust all locally functional
polygroupoids.

\begin{thm}[Group coset frame theorem, \cite{andgiv1, gaapal}]\label{lpg-thm} $\LPG$ is exactly the
class of structures associated to group coset frames, up to
isomorphisms. In more detail:
\begin{description}
%\item[(i)]
%For each group coset frame $(\G,\varphi,C)$, the associated
%structure
% $\Pp(\G,\varphi,C)$ is a locally functional polygroupoid.
\item[(i)]
$\Pp(\G,\varphi,C)$ is a locally functional polygroupoid if and only
if $(\G,\varphi,C)$ is a group coset frame.
\item[(ii)]
For each locally functional polygroupoid $\Pp$ there is a group
coset frame $(\G,\varphi,C)$ such that $\Pp$ is isomorphic to
$\Pp(\G,\varphi,C)$.
 \end{description}
 \end{thm}

\noindent{\bf Proof.} Proving part (i) of the theorem in relation
algebraic context is the main subject of \cite{andgiv1} and part
(ii) in relation algebraic form is the main representation theorem
of \cite{gaapal}. After showing how the present theorem follows from
the results in these papers, we sketch a direct proof of (ii) in
terms of polygroupoids.

 To be able to apply the above relation algebraic theorems,
 we begin by showing that $\Pp=\Pp(\G,\varphi,C)$ is isomorphic to the structure of
 atoms of the algebra $\Cc=\Cc(\G,\varphi,C)$ defined in
 \cite{andgiv1}, below Definition 3.2.

 The definition of $\Cc$ is as follows. First, we make the universes
 of the groups $\Gg_x$ disjoint. Let $(x,y)\in\E$ and let $H$ be a coset of $H_{xy}$ in $\Gg_x$.
 Define (see \cite[Definition 2.1]{andgiv1})
 \[ R_{xy,H} = \bigcup\{ K\times\varphi_{xy}[K\scir H] : K\in
 G_x\slash H_{xy} \} .\]
 The atoms of $\Cc$ are the binary relations $R_{xy,H}$ and
 multiplication $\otimes$ of $\Cc$ between atoms is (see
 \cite[Definition 3.1]{andgiv1})
 \[ R_{xy,H}\otimes R_{yz,K} = \bigcup\{ R_{xz,Q} : Q\subseteq\varphi_{yx}[\varphi_{xy}(H)\scir K]\scir C_{xyz}\},\]
 and $R_{xy,H}\otimes R_{wz,K}=\emptyset$ when $w\ne y$.
 The identity atoms are $R_{xx,H}$ where $H=H_{xx}$, by
 \cite[Theorem 2.4]{andgiv1}.
The inverse of $R_{xy,H}$ is $R_{yx,K}$ where
 $K=\varphi_{xy}(H^{-1})$, see \cite[Theorem 2.5]{andgiv1}.
 Therefore, the function assigning $R_{xy,H}$ to $(x,H,y)\in P$ is
 the desired isomorphism between $\Pp$ and the structure of atoms of $\Cc$.

Theorem 4.2 in \cite{andgiv1} states that $\Cc$ is a complete atomic
measurable relation algebra when $(\G,\varphi,C)$ is a group coset
frame. Comer \cite[Theorem 4.2]{comer1} states that the structure of
atoms of a complete atomic relation algebra is a polygroupoid. Thus
$\Pp\in\PG$. It is locally functional because $\Cc$ is measurable,
but this can be seen easier directly, as follows. The loops among
the elements of $\Pp$ are of form $(x,g,x)$ where $x\in I$ and $g$
is a coset of $H_{xx}$. When $(\G,\varphi,C)$ is a group coset
frame, we have that $H_{xx}=\{ e_x\}$, $\varphi_{xx}$ is the
identity map on cosets of $H_{xx}$ and $C_{xxx}=H_{xx}$, so
$\varphi_{xx}(\varphi_{xx}(g)\scir h)\scir C_{xxx}=g\scir h$ which
is a single coset of $H_{xx}$ and therefore $(x,g,x)(x,h,x)=\{
(x,k,x) : k\subseteq\varphi_{xx}(\varphi_{xx}(g)\scir h)\scir
C_{xxx}\}$ consists of a single atom. Also,
$(x,g,x)(y,h,h)=\emptyset$ when $y\ne x$. Thus $\Pp$ is locally
functional. This proves one direction of (i). The other direction,
namely that $(\G,\varphi,C)$ is a group frame when $\Pp$ is a
locally functional polygroupoid follows from Theorems 2.4, 2.5, 2.8,
3.4, 3.5, 3.7, 3.8 in \cite{andgiv1}.

To prove (ii) of the theorem, let us assume that $\Pp$ is a locally
functional polygroupoid. Let $\Cc$ be its natural algebra of
complexes, defined in \cite{comer1} for example. Then $\Cc$ is a
complete atomic relation algebra, by \cite[Theorem 4.2]{comer1}.
That $\Pp$ is locally functional translates immediately to the
statement that $\Cc$ is measurable (for the definition of a
measurable relation algebra see \cite{ga02}). Now, the main
representation theorem of \cite{gaapal}, Theorem 7.4, states that
there is a group coset frame $(\G,\varphi,C)$ such that the
structure of atoms of $\Cc$ is $\Pp(\G,\varphi,C)$. This proves
(ii).

Next we sketch a direct proof of Theorem \ref{lpg-thm}(ii), in terms
of polygroupoids. This will also show how to extract the group coset
frame $(\G,\varphi,C)$ from the polygroupoid.

Assume that $\Pp=(M,\scir,I,{}^{-1})$ is a locally functional
polygroupoid.  For any $x,y\in I$ define $M_{xy} = \{ a\in M :
x\scir a\scir y=a\}$. Let $\E=\{ (x,y)\in I\times I :
M_{xy}\ne\emptyset\}$. It can be shown that $\E$ is an equivalence
relation on $I$.

Let $G_x=M_{xx}$. It can be shown that $\Gg_x=(G_x,\scir,\{
x\},{}^{-1})$ is a group for each $x\in I$, because $\Pp$ is locally
functional. This gives our system $\G=\langle\Gg_x : x\in I\rangle$
of groups.

Next we construct the isomorphism $\varphi_{xy}$ between quotients
of $\Gg_x$ and $\Gg_y$ for $(x,y)\in\E$. For any $a\in M_{xy}$ let
\[ H_{xy}(a) = \{ g\in G_x : g\scir a=a\}\quad\mbox{ and }\quad
K_{xy}(a) = \{ g\in G_y : a\scir g=a\}\] and let $\varphi_{xy}(a)$
take a coset $H$ of $H_{xy}$ to a coset $K$ of $K_{xy}$ just in case
$H\scir a= a\scir K$. It can be shown that $H_{xy}(a), K_{xy}(a)$
are normal subgroups of $\Gg_x$ and $\Gg_y$ respectively, and they
do not depend on the choice of $a$. Moreover,
$\varphi_{xy}(a):\Gg_x\slash H_{xy}(a)\to \Gg_y\slash K_{xy}(a)$ is
an isomorphism, but this isomorphism may be different for different
elements $a$. To select a coherent system of isomorphisms, let us
select a system $a=\langle a_{xy} : (x,y)\in\E\rangle$ of elements
such that $a_{xy}\in M_{xy}$, $a_{xx}=x$ and $a_{xy}^{-1}=a_{yx}$
for all $(x,y)\in\E$. Such a system of elements can be chosen, it is
called a semi-scaffold. Let us define $\varphi_{xy} =
\varphi_{xy}(a_{xy})$, and let $\varphi=\langle\varphi_{xy} :
(x,y)\in\E\rangle$. It can be proved that $\varphi$ satisfies frame
conditions (i), (ii).

Finally, we define $C_{xyz}$ so that frame conditions (iii), (iv)
hold, as follows.  Ideally, we would want that $a_{xz}\in
a_{xy}\scir a_{yz}$ also hold in the semi-scaffold (then $a$ is
called a scaffold), but this cannot always be achieved. We define
$C_{xyz}$ to make up for this, by letting $g\in G_x$ be such that
$g\scir a_{xz}\in a_{xy}\scir a_{yz}$ (one can prove that there is
such a $g$), and then $C_{xyz}=g\scir H_{xy}\scir H_{xz}$. Finally,
one can show that coset conditions (v)-(viii) hold because $\Pp$ is
a polygroupoid. \qed

The notion of the Cayley representation of a group can be
generalized to polygroupoids. Since we have multivalued composition
in polygroupoids, we will represent the elements by binary relations
in place of bijective functions (permutations) as in a
Cayley-representation. When $R$ is a binary relation, its inverse is
$R^{-1}=\{ (v,u) : (u,v)\in R\}$.

\begin{Def}[Representation of polygroupoid]\label{rep-def}
A \emph{representation} of a polygroupoid $\Mm=(M,\scir,I,{}^{-1})$
on a set $U$ is a function $\rep$ that assigns disjoint nonempty
binary relations on $U$ to elements of $M$, that is
\begin{description}
\item{}
$\rep(a)\subseteq U\times U$ is nonempty and $\rep(a)$ is disjoint
from $\rep(b)$ whenever $a\ne b$,
\end{description}
that satisfies the following three conditions for all $a,b\in M$.
\begin{description}
\item{}
$\rep(a)\mid\rep(b)=\bigcup\rep[a\scir b]$,
\item{}
$\rep(x)\subseteq\{(u,u) : u\in U\}$\quad when $x\in I$, and
\item{}
$\rep(a)^{-1}=\rep(a^{-1})$.
\end{description}
Let $\RPG$ denote the class of representable polygroupoids.\qed
\end{Def}

Groups and Brandt groupoids are all representable, and in fact,
representations can be put together from the Cayley representations
of the groups involved. Our original hope was to prove the same for
all locally functional polygroupoids. However, after proving several
representation theorems for special classes (see \cite{ga02, giv1,
agtams}), we found a nonrepresentable locally functional
polygroupoid. The main result of the present paper is to construct a
series of such nonrepresentable locally functional polygroupoids
(Theorem \ref{nonrep-thm}), and to use them for proving that there
are continuum many varieties of relation algebras of a special kind
(Theorem \ref{var-thm}).

To prove nonrepresentability of our polygroupoids, we will use
Theorem \ref{gf-thm} below which characterizes the representable
locally functional polygroupoids.

\begin{Def}[Scaffold, {\cite[Definition 2]{ga02}}]\label{scaffold-def}
Let $\Mm=(M,\scir, I, {}^{-1})$ be a polygroupoid. A scaffold in
$\Mm$ is a system $a=\langle a_{xy} : x,y\in I, x\scir M\scir
y\ne\emptyset\rangle$ of elements of $M$ that satisfies the
following, for all $x,y,z\in I$ such that $x\scir M\scir
y\ne\emptyset$, $y\scir M\scir z\ne\emptyset$.
\begin{description}
\item{}
$a_{xy}=x\scir a_{xy}\scir y$,\quad $a_{xx}\in I$,\quad
$a_{xy}=a_{yx}^{-1}$\quad and\quad
%\item{}
$a_{xz}\in a_{ xy}\scir a_{yz}$. \qed
\end{description}
\end{Def}
Let us call $(\G,\varphi)$ a \emph{group frame} if $(\G,\varphi,C)$
is a group coset frame where all the $C_{xyz}$s are the identity
cosets $H_{xy}\scir H_{xz}$. The following theorem says that
choosing the coset system $C$ in a nontrivial way is essential in
constructing nonrepresentable $\LPG$s.

\begin{thm}[Group frame theorem, \cite{giv1, gaapal}]\label{gf-thm}
Let $\Mm\in\LPG$. The following are equivalent.
\begin{description}
\item[(i)]
$\Mm$ is representable.
\item[(ii)]
$\Mm$ is isomorphic to $\Pp(\G,\varphi,C)$ for some group coset
frame $(\G,\varphi,C)$  where all the $C_{xyz}$s are the identity
cosets $H_{xy}\scir H_{xz}$.
\item[(iii)]
There is a scaffold in $\Mm$.
%\item{}
%Ide a concrete representation-t a Givi cikkbol.
\end{description}
\end{thm}

\noindent{\bf Proof.} Assume that $\Mm\in\LPG$ and let $\Cc$ be its
algebra of complexes. Then $\Cc$ is a complete atomic measurable
relation algebra, by the proof of Theorem \ref{lpg-thm}. Any
representation of $\Mm$ can be extended to a function defined on the
complexes, simply by defining $\rr(X)=\bigcup\{ \rr(x) : x\in X\}$.
It can be proved that $\rr$ so defined is a complete representation
of $\Cc$ as defined in the relation algebra literature, and any
complete representation of $\Cc$ is a representation of $\Mm$ when
restricted to the atoms, that is to say to the elements of $M$.
Thus, $\Cc$ is completely representable if and only if $\Mm$ is
representable. By the proof of Theorem \ref{lpg-thm}, we have that
$\Cc$ is a complete and atomic measurable relation algebra. Thus,
$\Cc$ is completely representable if and only if it has a scaffold,
by Corollary 7.7 and Theorem 7.8 of \cite{gaapal}. Finally, $\Cc$ is
completely representable if and only if its structure of atoms is
determined by a group frame, by Theorem 7.6 of \cite{gaapal} and
Theorem 4.2 of \cite{giv1}. \qed

\section{Construction of group frames}\label{frame-sec}

In this section, we construct a series of group frames, one for each
commutative group $\Ff$. A coset system $C$  will be added in the
next section. Our goal is to add a coset system $C$ such that the so
obtained group coset frame will determine a nonrepresentable
polygroupoid. This places restrictions on how we can define our
group frame. For example, any group coset frame $(\langle \Gg_x :
x\in I\rangle, \varphi, C)$ where $I$ has less than five elements
determines a representable polygroupoid, so we want to have at least
five nodes $x$ in the group system. When all the groups $\Gg_x$ are
the product of at most two finite cyclic groups, we cannot achieve
nonrepresentability, either. Therefore we will have $\Ffff$ at the
nodes of our group frames. Finally, we have to have ``three distinct
levels" in the construction, because if at each node $x$ we have
$H_{xy}=\{ e_x\}$, or we have $H_{xy}=H_{xz}$ for all distinct
$y,z,x\in I$, then nonrepresentability cannot be achieved. Therefore
we want to have that $H_{xy}\ne\{ e_x\}$  and $H_{xy}\scir H_{xz}\ne
H_{xy}$ for all distinct $x,y,z$. These ``guiding" theorems are
mentioned in \cite{ga02}, and their proofs are in \cite{gamanu}.

In this section, we use additive notation for groups instead of
multiplicative notation as in the rest of the paper. Let
$\Ff=\langle F,+\rangle$ be a commutative group with zero element
$0$. We denote the inverse of $g\in F$ by $-g$. We shall define a
particular group frame $(G,\varphi)$ with all the $\Gg_x$ being
$\Ffff$.

\subsection{The normal subgroups}
Define the following subsets of $F\times
F\times F$:
\begin{description}
\item{}
$L_0 = \{ (a,0,0) : a\in F\}$,
\item{}
$L_1 = \{ (0,a,0) : a\in F\}$,
\item{}
$L_2 = \{ (0,0,a) : a\in F\}$,
\item{}
$L_3 = \{ (a,a,a) : a\in F\}$.
\end{description}
The properties of these subsets we will use in the paper are
gathered in the following two lemmas. Below,
$[u,v]=\{ (u,v), (v,u)\}\cup\{ (w,w) : w\in U\}$ denotes the
\emph{transposition} of $u,v$ whenever a set $U$ is understood from
context.

%\tag{i} -t valahogy boldfacce tenni
\begin{lem} \label{resp-lem}
For all permutations $\pi$ of $\{ L_m : m<4\}$ there is an
automorphism $\alpha(\pi)$ of $\Ffff$ which takes the subsets $L_m$
to their $\pi$-values, i.e.,
\[\tag{i}\label{aut-eq} \alpha(\pi)[L_m]=\pi(L_m)\quad\mbox{ for all }m<4.\]
Further, a system of such automorphisms can be chosen such that they
respect some structure of the permutations. Namely, let $\Pi$ denote
the set of permutations of  $\mathcal{L}= \{ L_m : m<4\}$, then
there is a system $\langle\alpha(\pi) : \pi\in \Pi\rangle$ such that
besides the above \eqref{aut-eq} for all $\pi, \sigma\in\Pi$ we have
\begin{description}
\item[(ii)]
The identity permutation of $\{ L_m : m<4\}$ is taken to the
identity automorphism of $\Ffff$.
\item[(iii)]
Inverse is respected: $\alpha(\pi^{-1})=\alpha(\pi)^{-1}$.
\item[(iv)]
Composition is respected:
$\alpha(\pi)\mid\alpha(\sigma)=\alpha(\pi\mid\sigma)$.
\item[(v)]
Transpositions are respected in the sense that for distinct
$m,n,p<4$ we have
\[ \alpha([L_m,L_n])\ \mbox{ \ is the identity on }L_p.\]
\end{description}
\end{lem}

\noindent {\bf Proof.} Instead of permutations of $\{ L_i : i<4\}$
we will work with permutations of $4 = \{ i : i<4\}$ (by identifying
a $\pi\in\Pi$ with the permutation that takes $i$ to $j$ just in
case $\pi(L_i)$ is $L_j$).

All permutations of 4 can be obtained as compositions of the
following four transpositions: $[0,1], [0,2], [0,3], [1,2]$. Indeed,
all permutations can be obtained as a composition of transpositions,
the transpositions of 4 not listed above are $[1,3]$ and $[2,3]$.
However, $[1,3]=[0,1]\mid [0,3]\mid [0,1]$ and $[2,3]$ can be
obtained similarly.

Our plan is the following. First we define the automorphisms for the
above four transpositions and then we define the automorphism
belonging to an arbitrary permutation  $\pi$ as the composition of
these automorphisms according to an arbitrary decomposition of $\pi$
to transpositions. We will show that this definition is correct
because the resulting automorphism does not depend on which
decomposition of the permutation we choose.

For $i<4$ and $a\in F$ we define $g_i(a)\in L_i$ as follows:
\begin{description}
\item{} $g_0(a)=(a,0,0)$,
\item{} $g_1(a)=(0,a,0)$,
\item{} $g_2(a)=(0,0,a)$,
\item{} $g_3(a)=(-a,-a,-a)$.
\end{description}
Then $L_i=\{ g_i(a) : a\in F\}$. Let $\alpha$ be an automorphism of
$\Ffff$ and let $\pi$ be a permutation of 4. We say that
\emph{$\alpha$ is good for $\pi$} if for all $i<4$ and $a\in F$ we
have
\[\alpha(g_i(a)) = g_j(a)\qquad \mbox{ for }j=\pi(i). \]
For $i<j<3$ we define $\alpha([i,j])$ as the automorphism of $\Ffff$
that interchanges the $i$-th and $j$-th coordinate in an element of
$F\times F\times F$, and we define $\alpha([0,3])$ as below. Thus,
by writing $\alpha_{ij}$ in place of $\alpha([i,j])$ we define
\begin{description}
\item{} $\alpha_{01}(a,b,c) = (b,a,c)$,
\item{} $\alpha_{02}(a,b,c)=(c,b,a)$,
\item{} $\alpha_{12}(a,b,c)=(a,c,b)$ .
\item{} $\alpha_{03}(a,b,c) = (-a, b-a, c-a)$ .
\end{description}
We have to show that $\alpha_{03}$ is an automorphism of $\Ffff$.
Clearly, $\alpha_{03}$ is a bijection on $F\times F\times F$ because
$\Ff$ is a group: Assume that $(-a,b-a,c-a)=(-a',b'-a',c'-a')$. Then
$-a=-a'$ so $a=a'$. Thus $b=b'$ by $b-a=b'-a'$, and similarly
$c=c'$. To show that $\alpha_{03}$ is a homomorphism, we have to use
that $\Ff$ is commutative. %Indeed,
%$\alpha_{03}(a+d,b+e,c+f)=(-(a+d), b+e-(a+d), c+f-(a+d))$.
Since $\Ff$ is commutative, we have that $b+e-(a+d)=b-a+(e-d)$ and
$c+f-(a+d)=c-a+(f-d)$, so
$\alpha_{03}(a+d,b+e,c+f)=(-(a+d),b+e-(a+d),
c+f-(a+d))=(-a+-d,b-a+(e-d),c-a+(f-e))=(-a,b-a,c-a)+(-d,e-d,f-d)=\alpha_{03}(a,b,c)+\alpha_{03}(d,e,f)$.
Thus, $\alpha_{03}$ is indeed an automorphism of $\Ffff$.

We show that for all $i<j<4$ we have
\begin{description}\label{L1}
\item{(L1)}
$\alpha_{ij}$ is good for $[i,j]$.
\end{description}
Indeed,
\begin{description}\item{} $\alpha_{01}(g_0(a))=\alpha_{01}(a,0,0)=(0,a,0)=g_1(a)$ and
$1=[0,1](0)$. Also,
\item{} $\alpha_{01}(g_1(a))=\alpha_{01}(0,a,0)=(a,0,0)=g_0(a)$
and $0=[0,1](1)$. Finally,
\item{}
$\alpha_{01}(g_2(a))=(0,0,a)=(0,0,a)=g_2(a)$ and
\item{} $\alpha_{01}(g_3(a))=(-a,-a,-a)=(-a,-a,-a)=g_3(a)$, and
\item{} $2=[0,1](2)$,\qquad $3=[0,1](3)$.
\end{description}
The proofs for $\alpha_{02}$ and $\alpha_{12}$ are completely
analogous. The proof for $\alpha_{03}$ is as follows.
\begin{description}
\item{}$\alpha_{03}g_0(a)=\alpha_{03}(a,0,0)=(-a,0-a,0-a)=(-a,-a,-a)=g_3(a)$,
\item{}$\alpha_{03}g_1(a)=\alpha_{03}(0,a,0)=(-0,a-0,0-0)=(0,a,0)=g_1(a)$,
\item{}$\alpha_{03}g_2(a)=\alpha_{03}(0,0,a)=(-0,0-0,a-0)=(0,0,a)=g_2(a)$,
and
\item{}
$\alpha_{03}g_3(a)=\alpha_{03}(-a,-a,-a)=(-(-a),-a-(-a),-a-(-a))=(a,0,0)=g_3(a)$.
\end{description}
This finishes the proof of (L1). Next we prove
\begin{description}
\item{(L2)}
$\alpha\mid\beta$ is good for $\pi\mid\sigma$ whenever $\alpha$ and
$\beta$ are good for $\pi,\sigma$ respectively.
\end{description}
Indeed,
$(\alpha\mid\beta)(g_i(a))=\beta(\alpha(g_ia))=\beta(g_{\pi(i)}(a))
= g_{\sigma(\pi(i))}(a)$, so (L2) holds.
\begin{description}
\item{(L3)}
$\alpha=\beta$ whenever they are good for the same $\pi$.
\end{description}
Indeed, let $(a,b,c)\in F\times F\times F$ be arbitrary. Then
$(a,b,c)=(a,0,0)+(0,b,0)+(0,0,c)=g_0(a)+g_1(b)+g_2(c)$. Since
$\alpha,\beta$ agree on $g_0(a),g_1(b),g_2(c)$ and both are
homomorphisms, they agree on $(a,b,c)$, so (L3) holds.

We are ready to define $\alpha(\pi)$ for arbitrary permutations
$\pi$ of 4. Let us take two decompositions of $\pi$ into
transpositions, say
\[\pi\ \ =\ \ [i_1,j_1]\mid\dots\mid[i_p,j_p]\ \ =\ \ [m_1,n_1]\mid\dots\mid[m_t,n_t].\]
By (L1)-(L3) we have that
$\alpha([i_1,j_1])\mid\dots\mid\alpha[(i_p,j_p)]=\alpha([m_1,n_1])\mid\dots\mid\alpha[(m_t,n_t)]$,
so the following definition is correct.
\[\alpha(\pi)=\alpha([i_1,j_1])\mid\dots\mid\alpha[(i_p,j_p)].\]
By (L1),(L2) we have that $\alpha(\pi)$ is good for $\pi$. This
immediately implies that (i) and (v) of Lemma \ref{resp-lem} hold.
By the definition of $\alpha(\pi)$ we immediately get that
composition is respected, i.e., (iv) of Lemma \ref{resp-lem} holds.
Let $\gamma$ denote the identity permutation of 4 and let
$\mbox{id}$ denote the identity automorphism of $\Ffff$. Then both
$\mbox{id}$ and $\alpha(\gamma)$ are good for $\gamma$, hence they
are equal by (L3). This shows that identity is respected. This in
turn implies that inverse is also respected, since we have
\[\alpha(\pi)\mid\alpha(\pi^{-1})=\alpha(\gamma)=\mbox{id},\] therefore
$\alpha(\pi^{-1})$ is the inverse of the automorphism $\alpha(\pi)$.
Thus Lemma \ref{resp-lem}(ii),(iii) also hold.
The proof for Lemma \ref{resp-lem} is complete. \qed

\begin{lem} \label{L-lem} Assume $i,j,k<4$ are distinct.
\begin{description}
\item[(i)]
$L_i$ is a normal subgroup of $\Ffff$.
\item[(ii)]
$L_i+L_j$ intersected with $L_i+L_k$ is $L_i$.
\item[(iii)]
$L_i+L_j+L_k=F\times F\times F$ .
\item[(iv)]
Let $C,D$ be cosets of $L_i+L_j$ and $L_i+L_k$ respectively. Then
$C\cap D$ is a coset of $L_i$.
\item[(v)]
$\Ffff\slash L_i \cong \Fff$, \ \ and\ \ $\Ffff\slash(L_i+L_j) \cong
\Ff$.
\end{description}
\end{lem}
\noindent {\bf Proof.} It is enough to check Lemma \ref{L-lem} for
$(i,j,k)=(0,1,2)$, because after that we can use Lemma
\ref{resp-lem} and the fact that all the statements in the lemma are
preserved by automorphisms to infer that lemma \ref{L-lem} holds for
all distinct $i,j,k<4$.

Clearly, $L_0$ is closed under addition, so it is a subgroup of
$\Ffff$, and it is normal because $\Ff$ is commutative, so $\Ffff$
is also commutative.

We have $L_0+L_1=\{(a,b,0) : a,b\in F\}$ and $L_0+L_2=\{(a,0,c) :
a,c\in F\}$, their intersection is $L_0=\{ (a,0,0) : a\in F \}$.

Clearly, $L_0+L_1+L_2=\{(a,b,c) : a,b,c\in F\}=F\times F\times F$.

Let $C,D$ be cosets of $L_0+L_1$ and $L_0+L_2$ respectively. Then
$C=\{ (a,b,c) : a,b\in F\}$ for some $c\in F$ and $D=\{ (a,b,c) :
a,c\in F\}$ for some $b\in F$. Their intersection is $C\cap
D=\{(a,b,c) : a\in F\}$ which is the coset $L_0+\{(0,b,c)\}$ of
$L_0$.

Finally, the elements of $(F\times F\times F)\slash L_0$ are
$\{(a,b,c) : a\in F\}$ for some $b,c\in F$, and it is easy to see
that the function assigning $(b,c)$ to $\{(a,b,c) : a\in F\}$ is an
isomorphism from $\Ffff\slash L_0$ onto $\Fff$. The proof for
$\Ffff\slash(L_0+L_1)\cong\Ff$ is similar.

This proves Lemma \ref{L-lem}. \qed

\subsection{The system of quotient isomorphisms}
Let $I$ be any five-element set. Define \[ \Gg_x = \Ffff,\qquad
\mbox{ for all }x\in I ,\] and $\G=\langle \Gg_x : x\in I\rangle$.
Let us specify $H_{xy}$ for distinct $x,y\in I$ such that
\begin{equation*}\tag{1}\label{h-eq}
\{ H_{xy} : y\in I, y\ne x\} = \{ L_i : i<4\}\qquad\mbox{ for any
}x\in I.
\end{equation*}
This can be done because $I$ has more than five elements. Since $I$
has five elements, we have that $H_{xy}$ and $H_{xz}$ are distinct
for distinct $x,y,z$.

Let $x,y\in I$ be distinct. We define a function $\pi_{xy}$ on $\{
L_i : i<4\}$ as follows.
\begin{equation*}\tag{2}\label{pi-eq}
\pi_{xy}(H_{xy}) = H_{yx}  \quad\mbox{ and }\quad
\pi_{xy}(H_{xu})=H_{yu}\ \ \mbox{ for any }u\in I\setminus\{ x,y\}.
\end{equation*}

\begin{lem}\label{pi-lem} Let $x,y,z\in I$ be distinct.
\begin{description}
\item[(i)]
$\pi_{xy}$ is a permutation of $\{ L_i : i<4\}$.
\item[(ii)]
$\pi_{yx} = \pi_{xy}^{-1}$.
\item[(iii)]
$\pi_{xy}\mid \pi_{yz} = \pi_{xz}\mid [H_{zx},H_{zy}]$.
\end{description}
\end{lem}

\noindent{\bf Proof.} Let $x,y,u,v,w$ be a repetitionfree listing of
$I$. Then $\pi_{xy}$ takes \[H_{xy}, H_{xu}, H_{xv},
H_{xw}\quad\mbox{ to }\quad H_{yx}, H_{yu}, H_{yv}, H_{yw},\]
respectively by \eqref{pi-eq}. Both lists are repetition-free
listings of $\{ L_i : i<4\}$, by \eqref{h-eq}. This shows that
$\pi_{xy}$ is a permutation of $\{ L_i : i<4\}$.

To prove (ii), let $p\in I$ be distinct from $x,y$. Then
\begin{description}
\item{}
$\pi_{yx}\pi_{xy}(H_{xy})=\pi_{yx}(H_{yx})=H_{xy}$ and
\item{}
$\pi_{yx}\pi_{xy}(H_{xp})=\pi_{yx}(H_{yp})=H_{xp}$.
\end{description}
Thus $\pi_{xy}\mid \pi_{yx}$ is the identity mapping on
$\{ L_i : i<4\}$. Together with (i), this proves (ii).

To prove (iii), let $p\in I$ be distinct from $x,y,z$ and let
$\pi_1$ denote $[H_{zx},H_{zy}]$. Now
\begin{description}
\item{}
$\pi_{yz}\pi_{xy}(H_{xy})=\pi_{yz}(H_{yx})=H_{zx}$ and
$\pi_1\pi_{xz}(H_{xy})=\pi_1(H_{zy})=H_{zx}$.
\item{}
$\pi_{yz}\pi_{xy}(H_{xz})=\pi_{yz}(H_{yz})=H_{zy}$ and
$\pi_1\pi_{xz}(H_{xz})=\pi_1(H_{zx})=H_{zy}$.
\item{}
$\pi_{yz}\pi_{xy}(H_{xp})=\pi_{yz}(H_{yp})=H_{zp}$ and
$\pi_1\pi_{xz}(H_{xp})=\pi_1(H_{zp})=H_{zp}$.
\end{description}
This shows that (iii) holds. \qed
\bigskip

Let $\langle\alpha(\pi) : \pi\in P\rangle$ be a system of
automorphisms of $\Ffff$ that satisfies the properties stated in
Lemma \ref{resp-lem}. Define
\begin{equation*}\tag{3}\label{a-eq}
\alpha_{xy} =\alpha(\pi_{xy}).\end{equation*} Then $\alpha_{xy}$ is
an automorphism of $\Ffff$ that takes $H_{xy}$ to $H_{yx}$. Hence,
$\alpha_{xy}$ induces an isomorphism $\varphi_{xy}$ between
$\Ffff\slash H_{xy}$ and $\Ffff\slash H_{yx}$ by setting
\begin{equation*}\label{phi-eq}\tag{4} \varphi_{xy}(g\slash H_{xy}) =
\alpha_{xy}(g)\slash H_{yx} .\end{equation*} For any $x\in I$, let
us define $\varphi_{xx}$ to be the identity mapping on $\Gg_x$. By
this, we have defined the system $\varphi = \langle\varphi_{xy} :
x,y\in I\rangle$. Figure 2 illustrates the definition.

\begin{figure}\label{gphi-fig}
\centering \small \psfrag{xy}[b][b]{$\Fff$}
%\psfrag{yx}[b][b]{$G_{yx}$} \psfrag{a}[b][b]{$\alpha_{xy}$}
\psfrag{yx}[b][b]{$\Fff$} \psfrag{a}[b][b]{$\alpha_{xy}$}
\psfrag{yxz}[t][t]{$\Ff$} \psfrag{x}[r][r]{$\Ffff$}
%\psfrag{y}[l][l]{$G_y$} \psfrag{yz}[l][l]{$G_{yz}$}
\psfrag{y}[l][l]{$\Ffff$} \psfrag{yz}[l][l]{}
%\psfrag{zy}[l][l]{$G_{zy}$} \psfrag{z}[t][t]{$G_z$}
\psfrag{zy}[l][l]{} \psfrag{z}[t][t]{$\Ffff$}
\psfrag{f}[b][b]{$\varphi_{xy}$} \psfrag{n}[b][b]{$\varphi^z_{xy}$}
\psfrag{xyz}[t][t]{$\Ff$}
\includegraphics[keepaspectratio, width=0.8\textwidth]{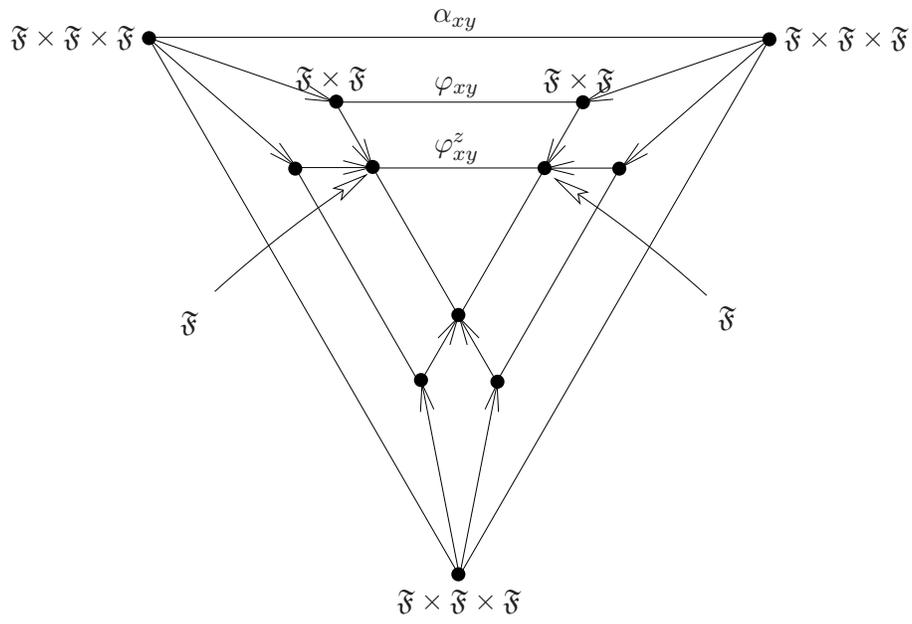}
\caption{The group frame $(\G,\varphi)$.}
\end{figure}

\begin{thm}\label{frame-thm}
$(\G,\varphi)$ is a group frame.
\end{thm}

\noindent{\bf Proof.} Throughout the proof, let $x,y,z\in I$ be
distinct  and $g\in F\times F\times F$ be arbitrary. $\varphi_{xx}$
is the identity by definition. We have $\pi_{yx}=\pi_{xy}^{-1}$ by
Lemma \ref{pi-lem}(ii), so $\alpha_{yx}=\alpha_{xy}^{-1}$ by
\eqref{a-eq} and Lemma \ref{resp-lem}(iii), so
$\varphi_{yx}=\varphi_{xy}^{-1}$ by \eqref{phi-eq}. Thus, conditions
(i) and (ii) of Definition \ref{frame-def} hold.

To check conditions (iii) and (iv), first we prove
\[\tag{5}\label{dupla-eq} \varphi_{xy}[g\slash H_{xz}\slash H_{xy}] =
\alpha_{xy}(g)\slash H_{yz}\slash H_{yx} .\]
Indeed,
\begin{description}
\item{}
$\varphi_{xy}[g\slash H_{xz}\slash H_{xy}] = \alpha_{xy}[g\slash
H_{xz}]\slash H_{yx} = \alpha_{xy}[\{ g\}+H_{xz}]\slash H_{xy} = $
\item{}
$(\{\alpha_{xy}(g)\} + \alpha_{xy}[H_{xz}])\slash H_{yx} =
(\{\alpha_{xy}(g)\} + H_{yz})\slash H_{yx} = \alpha_{xy}(g)\slash
H_{yz}\slash H_{yx}.$
\end{description}
Above, the first three equalities hold by \eqref{phi-eq}, $g\slash H
= \{ g\} + H$, and $\alpha_{xy}$ being a homomorphism, respectively.
The fourth equality holds because $\alpha_{xy}[H_{xz}]=H_{yz}$ by
\eqref{a-eq}, \eqref{pi-eq} and Lemma \ref{resp-lem}(i). The last
equality holds by $g\slash H = \{ g\} + H$.

Condition (iii) in Definition \ref{frame-def} is a special case of
\eqref{dupla-eq} when we take $g$ to be $0$, because $0\slash H=H$
for any subgroup $H$. Thus, condition (iii) holds.

To check condition (iv), first we prove
\[\tag{6}\label{phidef-eq} \varphi_{xy}^{z}(g\slash(H_{xy}+H_{xz})) =
\alpha_{xy}(g)\slash(H_{yx}+H_{yz})  .\]
Indeed, \begin{description}
\item{}
$\varphi_{xy}^{z}(g\slash (H_{xy}+H_{xz})) = \varphi_{xy}(g\slash
H_{xz})+ H_{yz} = \alpha_{xy}(g)\slash H_{yx} + H_{yz} =$
\item{}
$\alpha_{xy}(g)\slash (H_{yx}+H_{yz})$.
\end{description}
The first equation holds by the definition of $\varphi_{xy}^{z}$ in
Definition \ref{frame-def}, the second equality holds by
\eqref{phi-eq}, and the last equality holds by $\bigcup g\slash
H\slash H' = g\slash (H+H')$. This proves \eqref{phidef-eq}.

We will need that the automorphism associated to the transposition
of $L_i$ and $L_j$ for distinct $i,j<4$ is the identity modulo
$L_i+L_j$.
\[\tag{7}\label{id-eq}
\alpha([L_i,L_j])(g)\slash(L_i+L_j) = g\slash (L_i+L_j) .\] Indeed,
let $\alpha_1$ denote the isomorphism associated to the
transposition of $L_i$ and $L_j$ in $\{ L_i : i<4\}$, that is to
say, $\alpha_1=\alpha([L_i,L_j])$ and let $Q$ denote $L_i+L_j$. Then
$\alpha_1[L_i]=L_j$ and $\alpha_1[L_j]=L_i$ by Lemma
\ref{resp-lem}(i) and the definitions of $\alpha_1$ and the
transposition $[L_i,L_j]$. Therefore, $\alpha_1[Q]=Q$.
Let $L_p$ be distinct from $L_i$ and $L_j$. Then $g=n+f$ for some
$n\in Q$ and $f\in L_p$, by Lemma \ref{L-lem}(iii) since $L_i$,
$L_j$ and $L_p$ are distinct. Also, $\alpha_1(f)=f$ by Lemma
\ref{resp-lem}(v). Hence
\begin{description}\item{} $\alpha_1(g)\slash Q= \alpha_1(n+f)\slash
Q= \alpha_1(n)\slash Q+\alpha_1(f)\slash Q=$ \item{} $n\slash
Q+f\slash Q = (n+f)\slash Q = g\slash Q$.
\end{description}
The first and last equalities hold by $g=n+f$, the second and fourth
ones hold since $\alpha_1$ is a homomorphism and the third equality
holds because we have $\alpha_1(n)\slash Q=n\slash Q$ by
$\alpha_1[Q]=Q$ and $n\in Q$. This proves \eqref{id-eq}.

We are ready to check condition (iv).
\begin{description}
\item{}
$\varphi_{yz}^{x}\varphi_{xy}^{z}(g\slash (H_{xy}+H_{xz}))=$
\item{}
$\varphi_{yz}^{x}(\alpha_{xy}(g)\slash(H_{yx}+H_{yz})) =
%\item{}
\alpha_{yz}\alpha_{xy}(g)\slash(H_{zx}+H_{zy}) =$
\item{}
$\alpha([H_{zx},H_{zy}])\alpha_{xz}(g)\slash(H_{zx}+H_{zy}) =
\alpha_{xz}(g)\slash(H_{zx}+H_{zy}) =$
\item{}
$\varphi_{xz}^{y}(g\slash H_{xy}+H_{xz})$.
\end{description}
The first two and the last equalities hold by \eqref{phidef-eq}. The
third equality holds by Lemma \ref{resp-lem}(iii) and Lemma
\ref{pi-lem}(iii). The fourth equality holds by \eqref{id-eq}.
 This finishes checking condition
(iv) and the proof of Theorem \ref{frame-thm} is complete. \qed

\section{Shifting, proof of nonrepresentability}\label{shift-sec}

In this section, we return to multiplicative as opposed to additive
notation for groups: we denote the binary operation of a group by
$\scir$, the neutral or identity element is denoted by $e$ and the
inverse of an element $g$ in a group is denoted by $g^{-1}$.

Throughout this section, let $\Ff$ be any \emph{non-trivial}
commutative group, assume that we are given a particular choice for
$\langle H_{xy} : x,y\in I\rangle$ that satisfies \eqref{h-eq}
stated in the previous section, and
let $\mathcal{F}=(\G,\varphi)$ denote the group frame we defined
from $\Ff$ in the previous section.
%In fact, $\Ff$ does not determine the group frame, it also depends on
%the particular distribution of the normal subgroups on the edges of $I\times I$.
%The resulting group frame will be different, but the resulting group relation algebras will
%all be isomorphic, see cite[new isomorphism paper].
%
Let us give names for the elements of the index set $I$, i.e., let
\[ I = \{ p,q,r,s,t\} .\]
We will define the coset system $C$ such that it will be non-trivial
only in the triangle $pqr$.

Let $f\in F\times F\times F$ be such that $f\notin H_{pq}\scir
H_{pr}$. There is such an element because $\Ffff\slash H_{pq}\scir
H_{pr}\cong \Ff$ by Lemma \ref{L-lem}(v), \eqref{h-eq} and because
we assume that $\Ff$ has more than one element. We now define three
elements associated to $f$ (we could say, three ``versions" of $f$):
\[ \tag{8}\label{f-def}f_p= f,\quad f_q=\alpha_{pq}(f),\quad \mbox{and}\quad f_r=
\alpha_{pr}(f) .\]
We call $(x,y,z)$  a \emph{positive} permutation of $(p,q,r)$ if
$(x,y,z)$ is one of $(p,q,r)$, $(q,r,p)$, $(r,p,q)$, otherwise we
call $(x,y,z)$ a \emph{negative} permutation of $(p,q,r)$.

Define $C=\langle C_{xyz} : x,y,z\in I\rangle$ as follows.
\begin{description}
\item{}
$C_{xyz} = f_x\slash(H_{xy}\scir H_{zz})$, when $(x,y,z)$ is a
positive permutation of $(p,q,r)$,
\item{}
$C_{xyz} = f_x^{-1}\slash(H_{xy}\scir H_{xz})$, when $(x,y,z)$ is a
negative permutation of $(p,q,r)$,
\item{}
$C_{xyz} = H_{xy}\scir H_{xz}$, when $(x,y,z)$ is not a permutation
of $(p,q,r)$.
\end{description}

\begin{thm}\label{gcf-thm}
$(\G,\varphi,C)$ is a  group coset frame.
\end{thm}

\noindent{\bf Proof.}
We have to show that $(\G,\varphi,C)$ satisfies the conditions
listed in Definition \ref{cosetframe-def}. Since $(\G,\varphi)$ is a
group frame by Theorem \ref{frame-thm}, frame conditions (i)-(iii)
hold. Also (iv) holds with $\tau(C_{xyz})$ being the identity
function on the cosets of $H_{xy}\scir H_{xz}$, by Theorem
\ref{frame-thm}. Since $\Ff$ is commutative by assumption, we have
that $\Ffff$ is also commutative, and thus all inner automorphisms
are the identity function in it. Thus, frame condition (iv) holds
for our particular choice of $C$, too.

We now check that $C$ satisfies the four coset conditions (v)-(viii)
listed in Definition \ref{cosetframe-def}. Our system of cosets has
the special property that only cosets $C_{xyz}$ with $\{ x,y,z\}=\{
p,q,r\}$ are not the identity cosets. This special type of coset
system is dealt with in \cite[Corollary 4.8]{andgiv1}, we will check
the four conditions listed in that corollary.

Let us order the elements of $I$ as $p<q<r<s<t$. The first
condition, condition (i) of \cite[Corollary 4.8]{andgiv1} asks that
$C_{pqr}^{-1}=C_{prq}$ and the analogues for $qrp, rpq$ hold.
Condition (i) of \cite{andgiv1} is satisfied because
\[ C_{pqr}^{-1}=(f_p\slash(H_{pq}\scir
H_{pr}))^{-1}=f_p^{-1}\slash(H_{pq}\scir H_{pr})=C_{prq} .\] The
cases $qrp, rpq$ are completely analogous. This proves (i) of
\cite{andgiv1}. Condition (ii) of \cite{andgiv1} asks that
$\varphi_{pq}[C_{pqr}]=C_{qrp}$. Indeed,
\begin{description}\item{}$ \varphi_{pq}[C_{pqr}]=
\varphi_{pq}[f_p\slash(H_{pq}\scir H_{pr})]=
\alpha_{pq}(f_p)\slash(H_{pq}\scir H_{pr})= f_q\slash(H_{pq}\scir
H_{pr})=$ \item{} $C_{qrp}$.
\end{description}
The first and last equalities hold by the definitions of $C_{pqr},
C_{rpq}$ respectively. The second equality holds by the definition
of $\varphi_{xy}$ in \eqref{phi-eq}, and the third equality holds by
the definition of $f_r$ in \eqref{f-def}. This shows that coset
condition (ii) of \cite{andgiv1} holds. Checking condition (iii) of
\cite{andgiv1} is similar:
\begin{description}\item{}$ \varphi_{pr}[C_{pqr}]=
\varphi_{pr}[f_p\slash(H_{pq}\scir H_{pr})]=
\alpha_{pr}(f_p)\slash(H_{pq}\scir H_{pr})= f_r\slash(H_{pq}\scir
H_{pr})=$ \item{} $C_{rpq}$.
\end{description}
Finally, condition (iv) of \cite{andgiv1} asks that
\begin{description}\item{} $C_{pqr}\subseteq H_{pq}\scir H_{pr}\scir
H_{ps}\cap H_{pq}\scir H_{pr}\scir H_{pt}.$\end{description} This
holds because $H_{pq}\scir H_{pr}\scir H_{px}=G_x=F\times F\times F$
for $x=s,t$ by Lemma \ref{L-lem}(iii) and \eqref{h-eq}. Thus, $C$ is
indeed a system of cosets that satisfies the four coset conditions
of \cite[Corollary 4.8]{andgiv1}. This completes the proof of
Theorem \ref{gcf-thm}. \qed

\begin{thm}\label{nonrep-thm}
The polygroupoid \ $\Pp(\G,\varphi,C)$  is not representable.
\end{thm}

\noindent{\bf Proof.} Theorem \ref{frame-thm} states that
$(\G,\varphi)$ is a group frame, let $\Mm=\Pp(\G,\varphi)$ be the
structure associated to $(\G,\varphi,H)$ with $H=\langle H_{xy}\scir
H_{xz} : x,y,z\in I\rangle$ the ``trivial" coset system, and let
$\Pp=\Pp(\G,\varphi,C)$. Both $\Pp$ and $\Mm$ are locally functional
polygroupoids, by Theorems \ref{lpg-thm}, \ref{frame-thm},
\ref{gcf-thm}. They have the same universe
\[ P=\{ (x,g,y) : x,y\in I\mbox{ and }g\in G_x\},\] they have the
same identity atoms and inverse operation, by Definition
\ref{poly-def}, only their multiplication operation differ, and that
differs only for elements $(x,g,y), (y,h,z)\in P$ with $\{
x,y,z\}=\{ p,q,r\}$. Let $\otimes$ denote multiplication of $\Pp$,
and let $*$ denote multiplication of $\Mm$.

Let us inspect these two multiplications. Let $(x,g,y), (y,h,z)\in
P$ and let us define
\[ g(xy)h = \varphi_{yx}(\varphi_{xy}(g)\scir h).\]
Note that $g(xy)h$ is a coset of $H_{xy}\scir H_{xz}$, because $g$
is a coset of $H_{xy}$, so $\varphi_{xy}(g)$ is a coset of $H_{yx}$,
thus $\varphi_{xy}(g)\scir h$ is a coset of $H_{yx}\scir H_{yz}$
since $h$ is a coset of $H_{yz}$. Therefore
$\varphi_{yx}(\varphi_{xy}(g)\scir h)$ is a coset of $H_{xy}\scir
H_{xz}$ by condition (iii) in Definition \ref{cosetframe-def}. Thus,
$g(xy)h$ is a coset of $H_{xy}\scir H_{xz}$, and so is then
$g(xy)h\scir C_{xyz}$.
For convenience, let us define
\[ (x,X,y) = \{ (x,g,y) : g\subseteq X\} \]
when $X$ is a coset of $H_{xy}\scir H_{xz}$ for some $z\in I$. With
this notation we have
\begin{description}
\item{}
$(x,g,y)*(y,h,z) = (x, g(xy)h, z)$,\quad and
\item{} $(x,g,y)\otimes(y,h,z) = (x, g(xy)h\scir C_{xyz}, z)$ .
\end{description}
By Theorem \ref{gf-thm} we have that $\Mm$ is representable. We call
$\Mm$ the \emph{representable, or unshifted, pair} of $\Pp$.

We prove that $\Pp$ is not representable by deriving a contradiction
from assuming that it is.
Assume that $\Pp$ is representable. By Theorem \ref{gf-thm} then it
has a scaffold, let us fix such a scaffold  $a$ for $\Pp$. Note that
the set of identity elements of $\Pp$ is $\{(x,e,x) : x\in I\}$, and
$\Pp$ is connected. Thus, $a=\langle a_{xy} : x,y\in I\rangle$ is a
system of elements of $P$ that satisfies for all $x,y,z\in I$
\begin{description}
\item{}
$a_{xy} = (x, g_{xy}, y)$  for some $g_{xy}$,
\item{}
$a_{xx} = (x,e,x)$,
\item{}
$a_{xy}=a_{yx}^{-1}$, \quad and
\item{}
$a_{xz}\in a_{xy}\otimes a_{yz}$.
\end{description}
We now check what this scaffold ``does" in the representable pair
$\Mm$ of $\Pp$.  Since $\Mm$ and $\Pp$ have the same universe,
\[\tag{9}\label{atom-eq} a_{xy}=(x,g_{xy},y)\quad\mbox{ is an element of }\Mm,\mbox{ too.} \]
Since composition $*$ of $\Mm$ coincides with composition $\otimes$
of $\Pp$ on elements not in the ``triangle $pqr$", we have
\[\tag{10}\label{compa-eq} a_{st}\in a_{sx}*a_{xt}\quad\mbox{ for
}x\in\{ p,q,r\},\] and also
\[\tag{11}\label{compb-eq} a_{pq}\in a_{px}*a_{xq}\quad\mbox{ for
}x\in\{ s,t\},\] and similarly for any $yz$ in place of $pq$ if
$y,z\in\{ p,q,r\}$. Since converse in $\Mm$ is the same as in $\Pp$
we continue to have
\[\tag{12}\label{conv-eq} a_{xy}=a_{yx}^{-1}\quad\mbox{ in }\Mm. \]
From this point on we will use that $\Mm$ is representable. Let
$\rr$ be any representation of $\Mm$. Now, $\rr(a_{st})\ne\emptyset$
by the definition of a representation, so there are $u_s, u_t$ such
that
\[ \tag{13}\label{st-eq}(u_s,u_t)\in\rr(a_{st}). \]
By \eqref{compa-eq} we have $\rr(a_{st})\subseteq
\rr(a_{sx})\mid\rr(a_{xt})$, thus by the definition of relation
composition $\mid$ there are $u_x$ for $x\in\{ p,q,r\}$ such that
\[ \tag{14}\label{pqr-eq} (u_s,u_x)\in \rr(a_{sx})\quad\mbox{ and }\quad
(u_x,u_t)\in\rr(a_{xt}) \] for all $x\in\{ p,q,r\}$. See Figure 3.

\begin{figure}\label{crown-fig}
\centering \small \psfrag{p}[r][r]{$u_p$} \psfrag{q}[b][b]{$u_q$}
\psfrag{r}[l][l]{$u_r$} \psfrag{s}[t][t]{$u_s$}
\psfrag{t}[t][t]{$u_t$}
\includegraphics[keepaspectratio, width=0.6\textwidth]{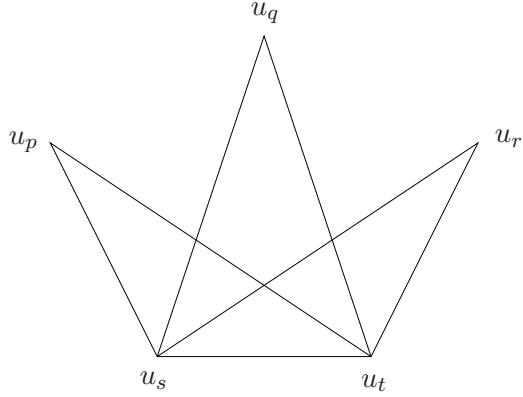}
\caption{Selection of the points in the representation of $\Mm$.}
\end{figure}

We now show that the structure of $\Mm$, that is the structure of
the group frame $(\G,\varphi)$, forces that $(u_x,u_y)\in
\rr(a_{xy})$ for $x,y\in\{ p,q,r\}$, and this will lead to a
contradiction.

We have $(u_p,u_s)\in \rr(a_{sp})^{-1}$ by $(u_s,u_p)\in
\rr(a_{sp})$. Also, $\rr(a_{sp})^{-1}=\rr(a_{sp}^{-1})=\rr(a_{ps})$
by the properties of $\rr$ and by \eqref{conv-eq} respectively. Now
$(u_s,u_q)\in \rr(a_{sq})$ implies that
$(u_p,u_q)\in\rr(a_{ps})\mid\rr(a_{sq})=\bigcup\rr[a_{ps}*a_{sq}]$.
Thus we got
\[ (u_p,u_q)\in\rr(b)\qquad\mbox{ for some } b\in a_{ps}*a_{sq}.\]
We get in a completely analogous way that
\[ (u_p,u_q)\in\rr(d)\qquad\mbox{ for some } d\in a_{pt}*a_{tq}.\]
Since representations of distinct elements are disjoint, we get that
$b=d$, and so
\[ \tag{15}\label{18}(u_p,u_q)\in\rr(b)\qquad\mbox{ for some } b\in (a_{ps}*a_{sq})\cap (a_{pt}*a_{tq}).\]
From \eqref{conv-eq} and \eqref{compb-eq} we also have that
\[\tag{16}\label{19} a_{pq}\in (a_{ps}*a_{sq})\cap(a_{pt}*a_{tq}).\]
We are going to show that $(a_{ps}*a_{sq})\cap(a_{pt}*a_{tq})$ is a
one-element set, this will imply $(u_p,u_q)\in \rr(a_{pq})$ by
\eqref{18}, \eqref{19}.

For all $x,y\in I$ we have that $a_{xy}=(x, g_{xy}, y)$ for some
coset $g_{xy}$ of $H_{xy}$. Thus $a_{ps}*a_{sq}=(p,D,q)$ for some
coset $D$ of $H_{pq}\scir H_{ps}$ and $a_{pt}*a_{tq}=(p,E,q)$ for
some coset $E$ of $H_{pq}\scir H_{pt}$. Since $H_{pq}, H_{ps},
H_{pt}$ are distinct elements of $\{ L_i : i<4\}$ by the
construction of the group frame $(\G,\varphi)$ and $D,E$ are cosets
of $H_{pq}\scir H_{ps}$ and $H_{pq}\scir H_{pt}$ respectively, Lemma
\ref{L-lem}(iv) implies that $D\cap E$ is a single coset of
$H_{pq}$. Therefore,  $b=(p, D\cap E,q)$ for any element
 $b$ of $(a_{ps}*a_{sq})\cap(a_{pt}*a_{tq})$. Thus indeed this set
 has a single element, and so $(u_p,u_q)\in\rr(a_{pq})$.

In a completely analogous way we get that $(u_p,u_r)\in a_{pr}$ and
$(u_r,u_q)\in a_{rq}$. Thus \[
(u_p,u_q)\in\rr(a_{pr})\mid\rr(a_{rq}) = \bigcup
\rr[a_{pr}*a_{rq}].\] Also $(u_p,u_q)\in \rr(a_{pq})$. By
disjointness of the representing relations we get that
\[\tag{17}\label{20} a_{pq}\in a_{pr}*a_{rq} .\] By the defining
conditions of our scaffold $a$ we also have
\[\tag{18}\label{21}
a_{pq}\in a_{pr}\otimes a_{rq} .\] However, we got $\otimes$ from
$*$ by shifting with a nonzero coset $K=f_p\slash (H_{pr}\scir
H_{pq})$. Thus we have that $a_{pr}*a_{pq}=(p,D,q)$ for some coset
$D$ of $H_{pr}\scir H_{pq}$, by the definition of $*$, and we have
that $a_{pr}\otimes a_{pq}=(p,D\scir K,q)$. Since $K$ is a nonzero
coset, we have that $D$ and $D\scir K$ are distinct cosets, so they
are disjoint. This in turn shows that \[ \mbox{$a_{pr}*a_{pq}$ is
disjoint from $a_{pr}\otimes a_{pq}$} \] which contradicts
\eqref{20}, \eqref{21}. We got a contradiction, so $\Pp$ cannot be
representable. \qed

\section{Continuum many subvarieties of coset relation
algebras}\label{var-sec}

In this section, we use the series of nonrepresentable locally
functional polygroupoids we constructed in sections \ref{frame-sec}
and \ref{shift-sec} to prove a theorem about relation algebras.

\emph{Relation algebras} are the subalgebras of complex algebras of
polygroupoids, \cite[Theorem 4.2]{comer1}. The definition of the
\emph{complex algebra}, or algebra of complexes, of a polygroupoid
is rather natural. The elements of the complex algebra are the
complexes, that is to say the subsets, of the structure, the
operations and relations of the structure extend naturally to
operations on these complexes, and additionally we take the Boolean
algebra structure of subsets of a set. Thus the complex algebra of a
polygroupoid-type algebra is a ``real" algebra, with totally defined
operations, see \cite[section 2]{comer1}. Complex algebras are one
of the key players of algebraic logic.

A relation algebra is \emph{representable} if it is a subalgebra of
the complex algebra of a representable polygroupoid, and it is a
\emph{coset relation algebra} if it is the subalgebra of the complex
algebra of a locally functional polygroupoid, up to
isomorphism. %
Both the class $\RRA$ of representable relation algebras and the
class $\CRA$ of coset relation algebras are varieties, a technical
term for being axiomatized by sets of equations, see \cite[Theorem
16.27]{giv18a}, \cite[Theorem 3.1]{gajsl}. We have
$\RRA\subseteq\CRA$ because the natural relation algebra of all
subsets of an equivalence relation is also a coset relation algebra.

It is proved in \cite[Theorem 4.10]{gajsl} that no variety between
\RRA\ and \CRA\ can be finitely axiomatized, and moreover, no such
variety can be axiomatized by a set of universal formulas containing
only finitely many variables \cite[Theorem 4.11]{gajsl}. No such
variety has decidable equational theory, either, by \cite[Theorem
8.5.(viii)]{tg}, or \cite[Corollary 2.9]{agnmem}.
%majd rendesebb  hivatkozast
The following theorem states that there are many varieties to which
the above theorems apply.

\begin{thm}\label{var-thm} There are continuum many varieties of
$\CRA$ that contain $\RRA$.
\end{thm}

\noindent{\bf Proof.} We follow the proof that J\'onsson gave for
\cite[Theorem 7.8]{J82}. An algebra is called \emph{simple} if it
has no nontrivial congruences. J\'onsson shows that if $\langle
\Aa_n : n\in J\rangle$ is a series of simple finite relation
algebras such that no $\Aa_n$ can be embedded into $\Aa_m$ for
distinct $n,m$, then the varieties generated by $\{ \Aa_n : n\in
S\}$ are all distinct for the subsets $S$ of $J$. We will use this
idea. However, since we also need that our varieties contain $\RRA$,
we will use a bigger variety than what $\{ \Aa_n : n\in S\}$
generates.

Let $\langle\Ff_n : n\in J\rangle$ be an infinite system of finite
commutative nontrivial groups of distinct size, e.g., we can take
$\Ff_n$ to be the $n$-element cyclic group, and take $J$ to be the
set of natural numbers $n\ge 2$. Let us fix now $n\ge 2$. Let
$\mathcal{F}_n$ be a group frame that we constructed from
$\Ff_n\times\Ff_n\times\Ff_n$ in section \ref{frame-sec}, let
$\Pp_n$ denote the nonrepresentable polygroupoid that we constructed
in section \ref{shift-sec}, and let $\Cc_n$ denote the complex
algebra of $\Pp_n$. We are going to show that $\langle \Cc_n : n\in
J\rangle$ is a system of simple finite nonrepresentable coset
relation algebras.

Since $\Ff_n$ is finite, by the construction we have that $\Pp_n$ is
also finite, so $\Cc_n$ is also finite. Since $\Pp_n$ is a locally
finite polygroupoid, we have that $\Cc_n$ is a coset relation
algebra. Next we show that $\Cc_n$ is nonrepresentable because
$\Pp_n$ is not representable and finite: Assume in contrary that
$\Cc_n$ is representable, say $\rr$ assigns binary relations to
elements of $\Cc_n$. Take the restriction of this representation to
the atoms of $\Cc_n$. Then we get a representation for the
polygroupoid $\Pp_n$ because $\rr(a)\mid\rr(b)=\bigcup\rr[a\scir b]$
by $\Cc_n$ being finite. Since $\Pp_n$ is not representable by
Theorem \ref{nonrep-thm}, we get that $\Cc_n$ is not representable.
Finally, $\Cc_n$ is simple because $\Pp_n$ is connected, see
\cite[Theorem 4.5]{comer1}.
\[\tag{1}\label{cc} \mbox{$\Cc_n$ is a finite, simple, nonrepresentable coset relation
algebra.}\] For any set $S\subseteq J$, let $\V(S)$ denote the class
of algebras in $\CRA$ into which none of $\Cc_n$ with $n\in S$ can
be embedded:
\[ \tag{2}\label{vs-eq}\V(S) = \{ \Aa\in\CRA : \Cc_n\mbox{ cannot be embedded into }\Aa,
\mbox{ for all }n\in S\} .\] By \cite[Theorem 7.1]{J82}, every
finite simple relation algebra is splitting in the class $\RA$ of
all relation algebras, thus the biggest variety of $\RA$ not
containing $\Cc_n$ is the class of all relation algebras into which
$\Cc_n$ cannot be embedded. This is called the conjugate variety of
$\Cc_n$, let us denote it by $\V^-(\Cc_n)$. By $\CRA\subseteq\RA$
then we have that
\[ \V(S) = \bigcap\{ \V^-(\Cc_n) : n\in S\}\cap\CRA .\]
This shows that $\V(S)$ is a variety, since it is the intersection
of varieties.

All the $\Cc_n$ are nonrepresentable by \eqref{cc}, therefore none
of them can be embedded into a representable relation algebra (since
$\RRA$ is closed under taking subalgebras). Also,
$\RRA\subseteq\CRA$, hence we get
\[ \RRA\subseteq \V(S)\subseteq\CRA.\]
It remains to show that all the $\V(S)$ are distinct. For this it is
enough to show, that \[\tag{3}\label{d-eq} \mbox{$\Cc_n$ cannot be
embedded into $\Cc_m$ if $n,m$ are distinct.}\] Indeed, assume that
\eqref{d-eq} holds and $S,Z$ are distinct subsets of $J$, say $n\in
S$ but $n\notin Z$. Then $\Cc_n\notin\V(S)$ since it can be embedded
into itself and $n\in S$. Also, $\Cc_n\in\V(Z)$ if \eqref{d-eq}
holds because $n\notin Z$ so no $\Cc_m$ with $m\in Z$ can be
embedded into $\Cc_n$. This shows that $\V(S)\ne\V(Z)$.

Next we show that \eqref{d-eq} holds.
Indeed, assume that $\Cc_n$ is embedded into $\Cc_m$, say by a
one-to-one homomorphism $f:\Cc_n\to\Cc_m$. We show that $n=m$. There
are five subidentity atoms both in $\Cc_n$ and in $\Cc_m$.
Consequently, the monomorphism $f$ has to take subidentity atoms to
subidentity atoms, because of the following. Since $f$ is a
homomorphism, it has to take the five subidentity atoms of $\Cc_n$
to five disjoint subidentity elements of $\Cc_m$, and since $f$ is a
monomorphism, it has to take an atom to a nonzero element. Since
$\Cc_m$ has only five subidentity atoms, this can be only if $f$
takes subidentity atoms to subidentity atoms.

Assume that $f(x)=y$ for a subidentity atom $x$ of $\Cc_n$. Then $y$
is a subidentity atom of $\Cc_m$, we have just showed this, and
$f(x\scir 1\scir x)=y\scir 1\scir y$. Now, the group
$\Ff_n\times\Ff_n\times\Ff_n$ can be recovered as the set of atoms
below $x\scir 1\scir x$. In particular, there are $n^3$ atoms below
$x\scir 1\scir x$ such that their sum is $x\scir 1\scir x$.
Similarly, there are $m^3$ atoms below $y\scir 1\scir y$ in $\Cc_m$.

Now, the atoms $g$ below $x\scir 1\scir x$, are special in the sense
that $g^{-1}\scir g=x$ is true for all of them (they are functional
in relation algebraic terminology). Thus, the same has to be true
for their images $f(g)$. It can be checked, that the functional
elements below $y\scir 1\scir y$ in $\Cc_m$ are exactly the atoms
below $y\scir 1\scir y$ (they are functional as in $\Cc_n$, but the
sum of more than one such atoms is never functional).
Thus, $f$ has to take the $n^3$ atoms below $x\scir 1\scir x$ to
$m^3$ atoms below $y\scir 1\scir y$ such that the sum of the images
is $y\scir 1\scir y$. This is possible only if in $\Cc_m$ also
$y\scir 1\scir y$ is the sum of $n^3$ atoms. Since $\Ff_n$ and
$\Ff_m$ have different finite cardinalities for $n\ne m$, this is
possible only if $n=m$. This proves \eqref{d-eq}, and with this the
proof of the theorem is also complete.\qed

\section*{Acknowledgement} We thank Robin Hirsch for calling our
attention to the connection between Brandt groupoids and group
relation algebras.

\bigskip\bigskip\bigskip

\noindent Alfr\'ed R\'enyi Institute of Mathematics,\\
Hungarian Academy of Sciences\\
Budapest, Re\'altanoda st.\ 13-15, H-1053 Hungary\\
andreka.hajnal@renyi.mta.hu, nemeti.istvan@renyi.mta.hu

\end{document}

%% file: article5a-def.tex
\newcommand{\comment}[1]{}

%\DeclareMathperator{\do}{do}
%\DeclareMathperator{\rg}{rg}

\newcommand{\dm}[1]{\mathrm{dm(#1)}}
\newcommand{\rg}[1]{\mathrm{rg(#1)}}

\newcommand{\conv}{{}^{\smallsmile}}
\newcommand{\inv}{{}^{-1}}
\newcommand{\rp}{\mathrel{\vert}}
\newcommand{\Pow}{\mathcal{P}}
\newcommand{\E}{\mathcal{E}}

\newcommand{\id}{\mbox{\sf Id}}
\newcommand{\di}{\mbox{\sf Di}}
\newcommand{\Z}{\mathcal{S}}
\newcommand{\HH}{\mathcal{H}}
\newcommand{\rep}{\mbox{\sf rep}}

\newcommand{\de}{\mbox{\ :=\ }}
\newcommand{\deiff}{\mbox{$\quad : \Leftrightarrow\quad$}}
\newcommand{\Aa}{\mbox{$\mathfrak A$}}
\newcommand{\Bb}{\mbox{$\mathfrak B$}}
\newcommand{\Cc}{\mbox{$\mathfrak C$}}
\newcommand{\Mm}{\mbox{$\mathfrak M$}}
\newcommand{\qed}{\hfill\mbox{$\Box$}\bigskip}
\newcommand{\rr}{\mbox{\sf rep}}

\newcommand{\Ff}{\mbox{$\mathfrak F$}}
\newcommand{\Fff}{\mbox{$\Ff\times\Ff$}}
\newcommand{\Ffff}{\mbox{$\Ff\times \Ff\times \Ff$}}
\newcommand{\Gg}{\mbox{$\mathfrak G$}}
\newcommand{\Pp}{\mbox{$\mathfrak P$}}
\newcommand{\Rr}{\mbox{$\mathfrak R$}}

\newcommand{\G}{\mbox{$\mathcal{G}$}}

\newcommand{\RA}{\mbox{\sf RA}}
\newcommand{\RRA}{\mbox{\sf RRA}}
\newcommand{\CRA}{\mbox{\sf CRA}}
\newcommand{\V}{\mbox{\sf V}}
\newcommand{\W}{\mbox{\sf W}}
\newcommand{\PG}{\mbox{\sf PG}}
\newcommand{\LPG}{\mbox{\sf LPG}}
\newcommand{\RPG}{\mbox{\sf RPG}}
\newcommand{\K}{\mbox{\sf K}}

\newcommand{\Cm}{\mbox{\bf Cm}}
\newcommand{\Sub}{\mbox{\bf S}}
\newcommand{\Hom}{\mbox{\bf H}}
\newcommand{\Pro}{\mbox{\bf P}}

\newcommand\scir{\raise2pt\hbox{$\,\scriptscriptstyle\circ\,$}}

\newtheorem{thm}{Theorem}[section]
\newtheorem{lem}{Lemma}[section]
\newtheorem{Def}{Definition}[section]
\newtheorem{rem}{Remark}[section]
\newtheorem{prb}{Problem}[section]
\newtheorem{cor}{Corollary}[section]

\newcommand{\la}{\lambda}